\documentclass[10pt, a4paper]{article}
\usepackage{graphicx}
\usepackage{color}
\usepackage{amsmath}
\usepackage{amssymb}
\usepackage{amscd}
\usepackage{bbm}
\usepackage{tikz}
\usepackage{framed}
\usetikzlibrary{positioning,arrows.meta,bending,calc}
\usepackage{hyperref}
\hypersetup{
    bookmarks=true,         
    colorlinks=true,       
    linkcolor=red,          
    citecolor=green,        
    filecolor=magenta,      
    urlcolor=cyan           
}
\usepackage[most]{tcolorbox}
\usepackage[utf8]{inputenc}
\usepackage{amsthm}
\usepackage{bbm} 
\usepackage{glossaries}
\usepackage[linesnumbered,lined,boxed,commentsnumbered]{algorithm2e}

\usepackage{marginnote}

\newtheorem{Theorem}{Theorem}

\newtheorem{Lemma}{Lemma}
\newtheorem{Proposition}{Proposition}
\newtheorem{Corollary}{Corollary}
\newtheorem{Remark}{Remark}

\newcounter{BigConst}                     
\newcommand{\nC}{                   
    \refstepcounter{BigConst}             
    \ensuremath{C_{\theBigConst}}
    }
\newcommand{\oC}[1]{\ensuremath{C_{\ref*{#1}}}}  

\newcounter{SmallConst}                     
\newcommand{\nc}{                   
    \refstepcounter{SmallConst}             
    \ensuremath{c_{\theSmallConst}}
    }
\newcommand{\oc}[1]{\ensuremath{c_{\ref*{#1}}}}  

\newcounter{gamma}                     

\newcounter{kappa}                     

\newcounter{delta}                     

\newcounter{theta}                     
\newcommand{\ntheta}{                   
    \refstepcounter{theta}             
    \ensuremath{\theta_{\thetheta}}
    }
\newcommand{\otheta}[1]{\ensuremath{\theta_{\ref*{#1}}}}  

\newcounter{L}                     

\newcounter{eps}                     

\newcommand{\norm}[1]{\left\|#1\right\|}%
\newcommand{\vertiii}[1]{{\left\vert\kern-0.25ex\left\vert\kern-0.25ex\left\vert #1 
    \right\vert\kern-0.25ex\right\vert\kern-0.25ex\right\vert}}

\newcommand{\vertiiii}[1]{{\left\vert\kern-0.25ex\left\vert\kern-0.25ex\left\vert\kern-0.25ex\left\vert #1 
    \right\vert\kern-0.25ex\right\vert\kern-0.25ex\right\vert\kern-0.25ex\right\vert}}

\newcommand\eps{\epsilon}

\newcommand{\beginproof}{{\bf Proof. {\hspace{0.2cm}}}}
\def \endproof
{{\mbox{}\nolinebreak\hfill\rule{2mm}{2mm}\par\medbreak}}

\DeclareMathOperator{\Span}{span}

\DeclareMathOperator*{\diam}{diam}

\DeclareMathOperator*{\Range}{Range}

\def\ds1{\textrm{1\kern-0.25emI}} 

\newcommand{\1}{\ensuremath{\mathbbm{1}}}

\newcommand \cA{{\cal A}}

\newcommand \cF{{\cal F}}

\newcommand \cN{{\cal N}}

\newcommand \bE{{\mathbb E}}

\newcommand \bN{{\mathbb N}}

\newcommand \bP{{\mathbb P}}

\newcommand \bR{{\mathbb R}}

\newcommand \bX{{\mathbb X}}

\newcommand \fm{{\mathfrak m}}

\newcommand{\vb}{{\boldsymbol{b}}}

\newcommand{\vlambda}{{\boldsymbol{\lambda}}}

\newcommand{\vv}{{\boldsymbol{v}}}

\newcommand{\vzero}{{\boldsymbol{0}}}
\newcommand{\vx}{{\boldsymbol{x}}}

\newcommand{\mysymbol}[3]{%
\newglossaryentry{#1}{%
      name={\ensuremath{#2}},%
      text={\ensuremath{#2}},%
      description={#3},%
      sort={#1}%
    }%
\expandafter\newcommand\expandafter{\csname smb#1\endcsname}{\gls{#1}}%
\expandafter\newcommand\expandafter{\csname #1\endcsname}{\ensuremath{#2}}%
}

\usepackage[english]{babel}

\usepackage{enumitem}

\begin{document}
\title{Upper bounds for the \texorpdfstring{$L^q$}{Lq} empirical process via generic chaining}
\author{Zong Shang\thanks{Department of Statistics, CREST-ENSAE, Institut Polytechnique de Paris, Palaiseau, France. Email: \texttt{zong.shang@ensae.fr}}}
\date{\today}
\maketitle

\begin{abstract}
    Using the generic chaining method, we derive upper bounds for the $L^q$ process of sub-Gaussian classes when $1 \leq q \leq 2$, thereby resolving an open problem posed by Al-Ghattas, Chen, and Sanz-Alonso in \cite{al-ghattas_sharp_2025}. Combined with the results of \cite{al-ghattas_sharp_2025}, this yields upper bounds for the $L^q$ process for all $1 \leq q < \infty$. We also present corollaries of this result in the geometry of Banach spaces, including high-probability bounds on the $\ell_q$-norm diameter of random hyperplane sections of convex bodies—where the subspaces are not necessarily uniformly distributed on the Grassmannian manifold—and the restricted isomorphic property for $\ell_q$ norm.
\end{abstract}

\section{Introduction}\label{sec_introduction}
Let $(\Omega,\cA,\mu)$ be a probability space, and $X$ be a random vector distributed as $\mu$. Let $1\le q<\infty$ be a real number. Let $\cF\subset L^q(\mu)$ be a class of real-valued functions defined on $\Omega$. We suppose $\vzero\in\cF$. We further suppose that there exists a distance $d_{\psi_2}$ on $\cF$ such that $\cF$ has sub-Gaussian increment with respect to $d_{\psi_2}$, that is, there exists some absolute constant $\nC\label{C_subgaussian}>1$, such that for any $f,g\in\cF$, and any $u>0$, we have $\bP(|f-g|(X)\geq u)\leq \oC{C_subgaussian}\exp(-\frac{u^2}{d_{\psi_2}(f,g)^2})$. 
Let $N\in\bN_+$ be a positive integer, and $X_1,\cdots,X_N$ be independent copies of $X$.

In this note, we are concerned with the high-probability upper bound of the following $L^q$ empirical process
\begin{align}\label{eq:objecitve_empirical_process}
    \sup\left( \left| \frac{1}{N}\sum_{i=1}^N\left|f(X_i)\right|^q - \bE \left|f(X)\right|^q\right| :\, f\in\cF \right).
\end{align}
This type of process has been extensively studied; to name a few, \cite{guedon_lp-moments_2007} investigated the case $q \ge 2$ when $\cF$ possesses a modulus of convexity, assuming that $\cF$ is a class of linear functionals $f(\cdot) = f_\vv(\cdot) = \langle \vv, \cdot \rangle$. \cite{zhivotovskiy_dimension-free_2024} considered the case where $\cF$ is an ellipsoid and consists of linear functions; his proof relies on a PAC-Bayesian argument, and thus it remains unclear how to extend it to general function classes. Such $L^q$ processes are extremely common in Geometric Functional Analysis. For instance, the $L^q(\mu) \hookrightarrow \ell_N^q$ embedding problem (see \cite[Section 14.8]{talagrand_upper_2014} and the references therein) essentially seeks upper bounds for the $L^q$ process—but in that context, it is only assumed that $\cF \subset L^q(\mu)$, without requiring $\cF - \cF \subset L_{\psi_2}$. In mathematical statistics and compressed sensing, the most important cases are $q=2$ and $q=1$: for example, \cite{plan_dimension_2014} applied the case $q=1$ to the problem of random hyperplane tessellations, while the case $q=2$, \cite{mendelson_reconstruction_2007,mendelson_empirical_2010,dirksen_tail_2015,bednorz_concentration_2014,mendelson_upper_2016}, due to its connection with the $\ell_2$ norm, manifests as the restricted isomorphic property \cite{candes_decoding_2005}, which plays a fundamental role in compressed sensing and statistical learning theory. The work \cite{al-ghattas_sharp_2025} is the most closely related to ours—they established an upper bound for the $L^q$ process when $q \ge 2$, and \cite{abdalla_dimension-free_2025} later provided an alternative simpler proof. In \cite[Remark 2.4]{al-ghattas_sharp_2025}, the authors posed the following open problem. 

\begin{framed}
    \centering
    What is the upper bound of the $L^q$ process when $1 < q < 2$?
\end{framed} Our main result resolves this problem.
Our main tool for investigating \eqref{eq:objecitve_empirical_process} is generic chaining. To this end, we introduce some notation. Let $(T,d)$ be a metric space. For any $n\in\bN$, let $N_n = 2^{2^n}$. We say a sequence of finite sets $(T_n)_{n=1}^\infty$ is admissible, if $|T_0|=1$, $|T_n|\leq N_n$, $T_n\subset T_{n+1}$, and $\cup_{n=0}^\infty T_n$ is dense in $T$ (with respect to the topology generated by $d$). We define Talagrand's $\gamma_2$ functional, \cite{talagrand_upper_2021}, by
\begin{align}\label{eq:def_gamma_2}
    \gamma_2(T,d) = \inf\bigg(\sup\left( \sum_{n=0}^\infty 2^{\frac{n}{2}}d(\vv,T_n):\, \vv\in T \right):\, (T_n)_{n=0}^\infty\mbox{ is admissbile}\bigg).
\end{align}

Denote the expectation with respect to empirical measure by $P_N$, that is, $P_N : f \in L^1(\mu) \mapsto \frac{1}{N}\sum_{i=1}^N f(X_i)$, and the expectation with respect to population measure by $P$, that is, $P : f \in L^1(\mu) \mapsto \bE f(X)$. 
Following the standard notation in empirical process theory, we denote \eqref{eq:objecitve_empirical_process} by $\sup\big(|(P - P_N)|f|^q|:\, f \in \cF\big)$. For convenience, we write $\gamma_2(\cF) = \gamma_2(\cF, d_{\psi_2})$, where $d_{\psi_2}$ is the sub-Gaussian metric on $\cF$, and $\diam(\cF) = \diam(\cF, d_{\psi_2})$.
The main conclusion of this paper is the following theorem. Its proof may be found in Section~\ref{sec:proof_theo:main}.
\begin{Theorem}\label{theo:main}
    Assume that $\cF$ is a class of functions containing $\vzero$ and having sub-Gaussian increments. Let $X_1,\cdots,X_N$ be independent copies of $X$. Then there exists an absolute constant $\nC\label{C_final}$ depending only on $q$ such that for any $u\geq 1$, there holds with probability at least $1-\exp(-u)$,
    \begin{align*}
        &\sup\big(|(P - P_N)|f|^q|:\, f \in \cF\big)\leq \oC{C_final}\bigg(\frac{\gamma_2^q(\cF)}{N^{\min\{1,q/2\}}} + \diam(\cF)^{q-1}\frac{\gamma_2(\cF)}{\sqrt{N}} + \diam(\cF)^q\left(\sqrt{\frac{u}{N}} + \frac{u^{q/2}}{N^{\min\{1,q/2\}}}\right)\bigg).
    \end{align*}
\end{Theorem}
In Remark~\ref{remark:optimality} below, we will see that Theorem~\ref{theo:main} is sharp.

We observe that Theorem~\ref{theo:main} exhibits a phase transition at $q = 2$. This arises because the margin of the $L^q$ process, $|(P_N - P)|f|^q|$, viewed as a sub-Weibull random variable of order $2/q$, is log-convex when $q > 2$ but log-concave when $q \le 2$. Consequently, the tail behavior in its concentration inequality changes, which in turn affects both the complexity term and the deviation term. Furthermore, the complexity term consists of two components: $\diam^{q-1}(\cF)\gamma_2(\cF)/\sqrt{N}$ and $\gamma_2^q(\cF)/N^{\min\{1, q/2\}}$. The phase transition appears only in the latter term, $\gamma_2^q(\cF)/N^{\min\{1, q/2\}}$. This is because the margin of the $L^q$ process (as we shall see below) exhibits a mixed tail behavior—when $q \ge 1$, its tail is sub-Gaussian for small deviations and sub-Weibull for large deviations. The Gaussian-type part has a tail probability independent of $q$, and it corresponds to the initial segment of the generic chaining, which produces the term $\diam^{q-1}(\cF)\gamma_2(\cF)/\sqrt{N}$; thus, this part does not undergo a phase transition as $q$ varies. In contrast, the sub-Weibull-type part has a tail probability that depends on $q$ (see \eqref{eq:Bernstein_one_function_sub_Weibull} below), and it corresponds to the terminal segment of the generic chaining, giving rise to the complexity term $\gamma_2^q(\cF)/N^{\min\{1, q/2\}}$, which therefore exhibits a phase transition with respect to $q$. As for the deviation term, it depends on the initial stage of the chain (roughly speaking, $n \sim \log_2(\lceil u \rceil)$), which, in the context of generic chaining, is controlled by the wimpy variance of the process (to use the terminology of \cite[pp.~314]{boucheron_concentration_2013}), and thus is also affected by a similar phase transition in the concentration inequality of sub-Weibull random variables.

\section{Applications}
In this section, we present several applications of Theorem~\ref{theo:main}, including the Restricted Isomorphic Property and a Dvoretzky–Milman type theorem.

\subsection{Restricted Isomorphic Property}

In this section, we apply Theorem~\ref{theo:main} to prove the Restricted Isomorphic Property (RIP).
The RIP has played a crucial role in mathematical statistics and compressed sensing over the past two decades, to name a few, \cite{candes_decoding_2005,foucart_mathematical_2013,adamczak_restricted_2011,lecue_geometrical_2024}. Roughly speaking, it characterizes an isomorphic correspondence between two norms, reflecting the isomorphism between the two normed spaces defined by the true (population) measure and the empirical measure. For example, consider a family of linear functionals $f(\cdot) = \langle \cdot, \vv \rangle$ in $\bR^d$, and consider $X$ to be an isotropic random vector in $\bR^d$ in the sense that $\bE[X\otimes X]:\vv\in\bR^d\mapsto \bE[\langle X,\vv\rangle X]\in\bR^d$ is the identity operator in $\bR^d$. Then $\tfrac{1}{N}\sum_{i=1}^N |f(X_i)|^q = \tfrac{1}{N}\sum_{i=1}^N |\langle X_i, \vv \rangle|^q = \tfrac{1}{N}\|\bX \vv\|_q^q$, where $\bX = [X_1 | \cdots | X_N]^\top \in \bR^{N \times d}$ is the design matrix. This empirical norm (raised to the $q$-th power) concentrates around its expectation, that is, the population norm $\bE |f(X)|^q = \|\langle X, \vv \rangle\|_{L^q}^q$, uniformly over a subset of $\bR^d$, where we abbrevaite $L^q(\mu)$ by $L^q$. In mathematical statistics and compressed sensing, we expect the observed data to faithfully reflect the underlying population structure, meaning that the two normed spaces $N^{-1/q}\|\bX \cdot\|_q$ and $\|\langle X, \cdot \rangle\|_{L^q}$ are isomorphic. When such an isomorphism typically holds on the entire space $\bR^d$, it is referred to as the \emph{Isomorphic Property} (IP); However, it usually holds only on a cone in $\bR^d$, especially when $d>N$, and is thus called the \emph{Restricted Isomorphic Property} (RIP).

Let $K$ be a convex body in $\bR^d$, let $G$ be a standard Gaussian random vector in $\bR^d$. Let $\ell_*(K) = \bE\sup(\langle\vv,G\rangle:\,\vv\in K)$ be the Gaussian mean width of $K$. Let $\diam(K) = \max(\|\vv\|_2:\, \vv\in K)$ be the $\ell_2$ diameter of $K$.
Let $R>0$ be a real number and define $RS_{L^q}^d = \{\vv\in\bR^d:\, \|\langle X,\vv\rangle\|_{L^q} = R\}$. Define $\cF(R) = \{\langle \cdot,\vv\rangle:\, \vv\in K\cap RS_{L^q}^d\}$. Without loss of generality, we assume that this set is nonempty. The following corollary provides sufficient conditions for the RIP to hold on the cone $\mathrm{cone}(\cF(R))=\{\alpha\vv:\, \alpha\ge 1,\, \vv\in \cF(R)\}$.
\begin{Corollary}\label{coro:RIP}
    Assume that $X$ is a sub-Gaussian, isotropic random vector. There exist absolute constants $\ntheta\label{theta_fixed_point_1}$, $\nc\label{c_prob_RIP} = \oc{c_prob_RIP}(q)$ and $\nc\label{c_RIP} = \oc{c_RIP}(q)$, such that, provided $R$ satisfies the following inequality:
    \begin{align}\label{eq:fixed_point_equation}
        \begin{aligned}
            & \ell_*(K\cap RS_{L^q}^d)\leq \otheta{theta_fixed_point_1}RN^{\min\{\frac{1}{2}, \frac{1}{q}\}},
        \end{aligned}
    \end{align}then with probability at least $1-\exp(-\oc{c_prob_RIP}N^{\min\{1,\frac{2}{q}\}})$, for any $\vv\in \mathrm{cone}(\cF(R))$, there holds
    \begin{align}\label{eq:RIP}
        \oc{c_RIP}\|\langle X,\vv\rangle\|_{L^q}\leq \frac{1}{N^{\frac{1}{q}}}\norm{\bX\vv}_q \leq \oc{c_RIP}^{-1}\|\langle X,\vv\rangle\|_{L^q}.
    \end{align}
\end{Corollary}

\beginproof Since $X$ is a sub-Gaussian random vector, for any $1 \le q \le \infty$ there exists an absolute constant $C(q)$ depending only on $q$ such that for all $\vv \in \bR^d$: if $q \ge 2$, then $\|\langle X, \vv \rangle\|_{L^2} \le \|\langle X, \vv \rangle\|_{L^q} \le C(q)\|\langle X, \vv \rangle\|_{L^2}$; if $1 \le q \le 2$, then $\|\langle X, \vv \rangle\|_{L^2} \le C(q)\|\langle X, \vv \rangle\|_{L^q}$. Indeed, by the Paley–Zygmund inequality, there exists $0 < \varepsilon < 1$ depending only on $q$ such that $\bP(|\langle X, \vv \rangle| \ge \tfrac{1}{\sqrt{2}}\|\vv\|_2) > \varepsilon$. Hence, $\|\langle X, \vv \rangle\|_{L^q}^q \ge \bE\big[|\langle X, \vv \rangle|^q \1_{\{|\langle X, \vv \rangle| \ge \tfrac{1}{\sqrt{2}}\|\vv\|_2\}}\big] \ge (\tfrac{1}{\sqrt{2}})^q \|\vv\|_2^q \varepsilon = (\tfrac{1}{\sqrt{2}})^q \varepsilon \|\langle X, \vv \rangle\|_{L^2}^q,$
and therefore $\|\langle X, \vv \rangle\|_{L^2} \le \sqrt{2}\,\varepsilon^{-1/q}\|\langle X, \vv \rangle\|_{L^q}$. Therefore, $\diam(K \cap R S_{L^q}^d) \le C(q) R$.

By homogeneity, we apply Theorem~\ref{theo:main} to $\cF(R)$ with $u = \min\{(100\oC{C_final}C^q(q))^{-2} N,\; (100\oC{C_final}C^q(q))^{-2/q} N^{\min\{1,\, 2/q\}}\}$, and obtain that, with probability at least $1 - \exp(-u)$, the following holds:
\begin{align*}
    &\sup_{\vv\in K\cap RS_{L^q}^d}\left(\left| \frac{1}{N}\sum_{i=1}^N\left|\langle\vv,X_i\rangle\right|^q - R^q \right|\right) \leq \oC{C_final} \frac{\gamma_2^q(K\cap RS_{L^q}^d)}{N^{\min\{\frac{q}{2},1\}}} + \oC{C_final}C(q)^{q-1}R^{q-1}\frac{\gamma_2(K\cap RS_{L^q}^d)}{\sqrt{N}} + \frac{1}{50}R^q.
\end{align*}When $q\geq 2$, \eqref{eq:fixed_point_equation} implies that $\ell_*(K\cap RS_{L^q}^d)N^{-1/2}\leq \otheta{theta_fixed_point_1}RN^{\frac{1}{q}-\frac{1}{2}}\leq\otheta{theta_fixed_point_1}R$. When $1\leq q< 2$, we also have $\ell_*(K\cap RS_{L^q}^d)N^{-1/2}\leq \otheta{theta_fixed_point_1}R$.
By \eqref{eq:fixed_point_equation} together with Talagrand's majorizing measure theorem, there exists $\otheta{theta_fixed_point_1}>0$ such that, $\gamma_2^q(K\cap RS_{L^q}^d)\leq \frac{1}{50\oC{C_final}}R^qN^{\min\{1,\frac{q}{2}\}}$, and such that $\oC{C_final}C(q)^{q-1}\gamma_2(K\cap RS_{L^q}^d)/\sqrt{N}\leq \frac{1}{50}R$. Combining the above conditions and noting that $R = \|\langle X, \vv \rangle\|_{L^q}$ holds for any $\vv \in K \cap R S_{L^q}^d$, the proof is complete. Here, one may take $\oc{c_RIP} = 1/2^{1/q}$ and $\oc{c_prob_RIP} = \min\{(100\oC{C_final}C^q(q))^{-2/q},\, (100\oC{C_final}C^q(q))^{-2}\}$.
\endproof

In particular, if $\ell_*(S_{L^q}^d) \le \otheta{theta_fixed_point_1} N^{\min\{\frac{1}{2},\, \frac{1}{q}\}}$, then the smallest $R$ satisfying \eqref{eq:fixed_point_equation} is $0$. Consequently, with probability at least $1-\exp(-\oc{c_prob_RIP}N^{\min\{1,\frac{2}{q}\}})$, the RIP property \eqref{eq:RIP} holds uniformly for all $\vv \in \bR^d$, and in this case the cone $\mathrm{cone}(\cF(R))$ degenerates to the linear space $\bR^d$ itself, that is, there holds the IP.

\subsection{The \texorpdfstring{$\ell_p$}{ellp}-diameter of random sections of convex bodies}

Let $d > N$ be a positive integer. Let $X \in \bR^d$ be an isotropic sub-Gaussian random vector, and let $X_1, \ldots, X_N$ be independent copies of $X$. Define $\bX = [X_1 | \cdots | X_N]^\top \in \bR^{N \times d}$, which is a random matrix with independent rows, and let $E = \Range(\bX^\top) = \Span(X_1, \ldots, X_N)$. For any convex body $K$, let $\|\cdot\|_K$ denote the norm whose unit ball is $K$. Define the polar body $K^\circ = \{\vv \in \bR^d : \sup(\langle \vv, \vx \rangle : \vx \in K) \le 1\}$; then $\|\cdot\|_{K^\circ}$ is the dual norm of $\|\cdot\|_K$. A classical line of research in Banach space geometry concerns the study of random sections of convex bodies, including the Dvoretzky–Milman theorem, Milman's $M^*$ estimate, and Gluskin's theorem, etc., see, for example, the standard references \cite{pisier_volume_1989,artstein-avidan_asymptotic_2015,tomczak-jaegermann_banach-mazur_1989,aubrun_alice_2017,vershynin_high-dimensional_2018}. In this subsection, we apply Theorem~\ref{theo:main} to derive upper bounds on the $\ell_p$ norm of random sections. We have the following corollary.
\begin{Corollary}\label{coro:DM}
    With the notation introduced above, there exist absolute constants $\nc\label{c_DM_dimension}$, $\nc\label{c_lower_DM}$, and $\nC\label{C_upper_DM}$ that depend only on $p$ such that the following holds.
    \begin{enumerate}
        \item When $1<p\leq 2$, and $N\leq\oc{c_DM_dimension} \big(\frac{\ell_*(K)}{\diam(K,\ell_2)}\big)^{\frac{p}{p-1}}$, then  with probability at least $1-\exp\big(-N^{\frac{2(p-1)}{p}}\big)$,
        \begin{align*}
            \forall \vlambda\in\bR^N,\, \|\vlambda\|_p=1,\, \norm{\bX^\top\vlambda}_{K^\circ}\leq \oC{C_upper_DM}\ell_*(K).
        \end{align*}
        \item When $2<p\leq\infty$, and $N \leq \oc{c_DM_dimension} \big(\tfrac{\ell_*(K)}{\diam(K, \ell_2)}\big)^2$, then with probability at least $1-\exp(-N)$,
        \begin{align*}
            \forall \vlambda\in\bR^N,\, \|\vlambda\|_p=1,\, \norm{\bX^\top\vlambda}_{K^\circ}\leq \oC{C_upper_DM}\ell_*^{\frac{2(p-1)}{p}}(K)\diam(K)^{\frac{2-p}{p}}.
        \end{align*}
    \end{enumerate}
    Furthermore, if $2<p\leq\infty$, $X\sim\cN(\vzero,I_d)$ and $N\leq \oc{c_DM_dimension}\big(\tfrac{\ell_*(K)}{\diam(K, \ell_2)}\big)^{\frac{p+2}{p}}$, then there exists an absolute constant $\nc\label{c_devi_2}<1$ such that with probability at least $1-\exp(-\oc{c_devi_2}N)$,
    \begin{align*}
        \forall \vlambda\in\bR^N,\, \|\vlambda\|_p=1,\, \oc{c_lower_DM}\ell_*(K)\leq \|\bX^\top\vlambda\|_{K^\circ}\leq \oC{C_upper_DM}\ell_*^{\frac{2(p-1)}{p}}(K)\diam(K)^{\frac{2-p}{p}}.
    \end{align*}
\end{Corollary}
\beginproof Notice that $\sup(\|\bX^\top\vlambda\|_{K^\circ}:\, \|\vlambda\|_p=1) = \sup(\|\bX\vv\|_q:\, \vv\in K)$ where $q=\frac{p}{p-1}$. Moreover, since $X$ is isotropic and is sub-Gaussian, $(\bE|\langle \vv,X\rangle|^q)^{1/q}\lesssim_p (\bE|\langle\vv,X\rangle|^2)^{1/2}=\|\vv\|_2$. By Theorem~\ref{theo:main} applied to $\cF = \{\langle\cdot,\vv\rangle:\, \vv\in K\}$, $u = N^{\min\{\frac{2}{q},1\}}$, there holds
\begin{align*}
    &\sup\left(\norm{\bX\vv}_q:\, \vv\in K\right)^q  \leq (1+\oC{C_final})N\diam(K)^q + \oC{C_final}\bigg(N^{\max\{1-\frac{q}{2},0\}}\gamma_2^q(K) + \sqrt{N}\diam(K)^{q-1}\gamma_2(K)\bigg).
\end{align*}We discuss the two cases $q \ge 2$ and $1 \le q < 2$ separately.
\begin{enumerate}
    \item When $2\leq q< \infty$, that is, $1< p\leq 2$. By $N\lesssim (\gamma_2(K)/\diam(K))^q$, we have $\sqrt{N}\diam(K)^{q-1}\gamma_2(K)\lesssim \gamma_2^q(K)N^{\frac{2-q}{2q}}\lesssim \gamma_2^q(K)$.
    \item When $1\leq q<2$, that is, $2<p\leq \infty$, and $N\lesssim (\gamma_2(K)/\diam(K))^2$, we have $N\diam(K)^q,$ $N^{1-q/2}\gamma_2^q(K), $ and $\sqrt{N}\diam(K)^{q-1}\gamma_2(K) $ are smaller than $ \gamma_2^2(K)\diam(K)^{q-2}$.
\end{enumerate}Then, by Talagrand's majorizing measure theorem for Gaussian processes, all the upper bounds in Corollary~\ref{coro:DM} are established.

When $p \ge 2$ and $X$ is Gaussian, the lower bound follows from a standard net argument. Since $p \ge 2$, we have $\tfrac{p + 2}{p} \le 2$. Hence, the stronger Dvoretzky condition $N \le \oc{c_DM_dimension}\big(\tfrac{\ell_*(K)}{\diam(K, \ell_2)}\big)^{\frac{p+2}{p}}$ guarantees the validity of the upper bound in item~\emph{2}. Let $\eps = \tfrac{1}{4\oC{C_upper_DM}} \big(\tfrac{\diam(K)}{\ell_*(K)}\big)^{\frac{p-2}{p}}$. Let $S_p^N = \{\vlambda \in \bR^N : \|\vlambda\|_p = 1\}$, and let $V_\eps$ be an $\eps$-net of $S_p^N$, i.e., for any $\vlambda \in S_p^N$, there exists $\pi\vlambda \in V_\eps\subset S_p^N$ such that $\|\vlambda - \pi\vlambda\|_p \le \eps$. By \cite[Lemma~4.10]{pisier_volume_1989}, we have $|V_\eps| \le (1 + 2/\eps)^N$. For any $\pi\vlambda \in V_\eps$, note that $\bX^\top(\pi\vlambda) \sim \|\pi\vlambda\|_2 G$, so $\|\bX^\top(\pi\vlambda)\|_{K^\circ}$ has the same distribution as $\|\pi\vlambda\|_2 \|G\|_{K^\circ}$, which is almost surely greater than $\|\pi\vlambda\|_p \|G\|_{K^\circ} = \|G\|_{K^\circ}$ (since $p \ge 2$). By the Gaussian Lipschitz concentration inequality (see \cite[Theorem~4.7]{pisier_volume_1989}) and the union bound, there exist absolute constants $c,c'<1$ such that if $N \le c \eps (\ell_*(K)/\diam(K))^2 = \frac{c}{4\oC{C_upper_DM}}(\ell_*(K)/\diam(K))^{\frac{p+2}{p}}$, then with probability at least $1 - \exp(-c'\eps(\ell_*(K)/\diam(K))^{\frac{p+2}{p}})$, for all $\pi\vlambda \in V_\eps$, we have $\|\bX^\top(\pi\vlambda)\|_{K^\circ} \ge (1/2)\ell_*(K)$. From the upper bound in item~\emph{2}, we know that for all $\vlambda - \pi\vlambda$, one has
\[
\|\bX^\top(\vlambda - \pi\vlambda)\|_{K^\circ} \le \oC{C_upper_DM}\|\vlambda - \pi\vlambda\|_p\, \ell_*^{\frac{2(p-1)}{p}}(K)\, \diam(K)^{\frac{2-p}{p}}.
\]
Hence, we obtain that for any $\vlambda\in S_p^N$, there holds $\|\bX^\top(\vlambda - \pi\vlambda)\|_{K^\circ}\leq \frac{1}{4}\ell_*(K)$, and hence we may take $\oc{c_lower_DM}=\frac{1}{4}$.

\endproof

We refer to Corollary~\ref{coro:DM} as a Dvoretzky–Milman type theorem because it shows that, as long as the dimension of the random subspace $E = \Range(\bX^\top)$ does not exceed a certain critical value, the upper bound on the $\ell_p$-norm diameter of its intersection with the convex body $K^\circ$ remains stable with high probability—that is, it no longer depends on the dimension of the subspace $E$. Unlike the classical Dvoretzky–Milman theorem, however, there are no matching upper and lower bounds here, and thus the section is not necessarily sandwiched between two $\ell_p$ spheres of nearly equal radii. In addition, the notion of critical dimension in this context differs from that in the classical Dvoretzky–Milman theorem. Specifically, we observe that this critical dimension is given by $(\ell_*(K)/\diam(K))^{\min\{2,\, p/(p-1)\}}$, which exhibits a phase transition at $p=2$, corresponding to the Euclidean version of the Dvoretzky–Milman dimension.

\section*{Acknowledge}
The author thanks Radosław Adamczak, Jiaheng Chen, and Sjoerd Dirksen for their valuable suggestions. Part of this work was completed during a visit to the University of Warsaw, and the author is grateful to Radosław Adamczak for his hospitality and to the Erasmus+ PhD Mobility program for financial support.

\section{Proof of Theorem~\ref{theo:main}}\label{sec:proof_theo:main}

The case $q \ge 2$ has already been established in \cite{al-ghattas_sharp_2025} and in \cite{abdalla_dimension-free_2025}; here we only need to prove the case $1 \le q \le 2$. Therefore, in this section, we prove the following theorem.
\begin{Theorem}[$1\leq q\leq 2$]\label{theo:1<q<2}
    Grant the same conditions as in  Theorem~\ref{theo:main} but with $1\leq q\leq 2$. There exists an absolute constant $\oC{C_final}$ depending only on $q$ such that for any $x\geq 1$, with probability at least $1-\exp(-x)$ there holds
    \begin{align}\label{eq:objective_tail_prob}
        \begin{aligned}
            &\sup\left( \left| \frac{1}{N}\sum_{i=1}^N\left|f(X_i)\right|^q - \bE \left|f(X)\right|^q\right| :\, f\in\cF \right)  \\
            &\leq \oC{C_final}\Bigg(\frac{\gamma_2^q(\cF,d_{\psi_2})}{N^{q/2}} + \diam(\cF,d_{\psi_2})^{q-1}\frac{\gamma_2(\cF,d_{\psi_2})}{\sqrt{N}} +  \diam(\cF,d_{\psi_2})^q\left(\sqrt{\frac{x}{N}} + \left(\frac{x}{N}\right)^{\frac{q}{2}}\right)\Bigg).
        \end{aligned}
    \end{align}Moreover, for any $1\leq p<\infty$,
    \begin{align}\label{eq:objective_moments}
        \begin{aligned}
            &\norm{\sup_{f\in\cF}\bigg|(P_N-P)|f|^q\bigg|}_{L^p} = \left(\bE\left(\sup_{f\in\cF}\left( \big|(P_N-P)|f|^q\big| \right)\right)^p\right)^{\frac{1}{p}}\\
            &\leq \oC{C_final}\left( \frac{\gamma_2^q(\cF,d_{\psi_2})}{N^{q/2}} + \diam(\cF,d_{\psi_2})^{q-1}\frac{\gamma_2(\cF,d_{\psi_2})}{\sqrt{N}} +  \diam(\cF,d_{\psi_2})^q\left(\sqrt{\frac{p}{N}} + \left(\frac{p}{N}\right)^{q/2} \right) \right).
        \end{aligned}
    \end{align}
\end{Theorem}

\subsection{Notation and Preliminaries}\label{sec_premiminaries}
In this section, we introduce some background materials that will be used in the proof of Theorem~\ref{theo:main}.

\paragraph{Notation}Let $\psi_\alpha(x) := 2^{x^\alpha}-1$ where $\alpha>0$, that is, the Orlicz $\psi_\alpha$-function. Then the sub-Gaussian increment implies that for any $f,g\in\cF$, we have $\bE\psi_2(\frac{|f-g|(X)}{d(f,g)})\leq 1$ (for a well-chosen constant $\oC{C_subgaussian}$). We define the Orlicz $\psi_\alpha$ norm by $\|\cdot\|_{\psi_\alpha}:f\in \cF\mapsto \inf( C>0:\, \bE\psi_\alpha(\frac{|f(X)|}{C})\leq 1)$. 
We let $L_{\psi_\alpha}:= \{f :\, \|f\|_{\psi_\alpha} < \infty\}$ and since $\mu$ is a probability measure, we also refer to $f\in L_{\psi_\alpha}$ be a $\psi_\alpha$ random variable. Such random variables are sometimes also referred to as sub-Weibull random variables of order $\alpha$. We denote by $d_{\psi_\alpha}$ the metric induced by the norm $\|\cdot\|_{\psi_\alpha}$. 

\subsection{Generic Chaining} Generic chaining is a method that discretizes the index set of a stochastic process according to the tail probabilities of the process margins, thereby applying the union bound in a multiscale manner, \cite{talagrand_upper_2021}. 

Briefly, to construct an upper bound for a stochastic process via generic chaining, one uses an admissible sequence to witness the growth of the process. According to the target $p$-th moment, we choose a starting level $\ell = \lfloor \log_2(p) \rfloor$, after which the admissible sequence grows from time $\ell$ and gradually covers the entire index set (by assumption, $\cup_n T_n$ is dense in $\cF$). Consequently, any index of the process can be expanded along these times as a chain: for any $n > \ell$, define $\pi_n f$ to be any element in $\cF$ such that $d(f,T_n) = d(f,\pi_n f)$, then for every $n > \ell$ we have $f = (f - \pi_n f) + (\pi_n f - \pi_{n-1} f) + \cdots + (\pi_{\ell+1} f - \pi_\ell f) + \pi_\ell f$. Thus, we need to control the increment of the margin on each chain link, i.e., on $\pi_n f - \pi_{n-1} f$, and then sum these sufficiently small increments using the triangle inequality. In this procedure, since at each time $n$ we must apply a union bound over at most $|T_{n-1}|\,|T_n| \le N_{n+1} = 2^{2^{n+1}}$ pairs of random variables, which requires the tail probability of the increment of the margin on $\pi_n f - \pi_{n-1} f$ to be at least $2^{-2^{n+1}}$ in order to balance this metric complexity. Hence one typically takes $2\exp(-u 2^n)$ with $u \ge 2$ as a deviation parameter, see Lemma~\ref{lemma:union_bound} below. This necessitates a careful analysis of the tail behavior of the process increments on $\pi_n f - \pi_{n-1} f$ to guarantee the admissible $\exp(-u 2^n)$ tail bound.

In Theorem~\ref{theo:main}, we study the $q$-th power of a sub-Gaussian stochastic process. Therefore, we begin by examining its margin's concentration inequality when viewed as a (sub-)Weibull random variable.
\paragraph{Some facts about the \texorpdfstring{$\psi_\alpha$}{psi alpha} random variables.}We list some standard facts concerning $\psi_\alpha$ random variables. They can be found, for example, in \cite{sambale_notes_2023, kuchibhotla_moving_2022}.
If $\zeta \in L_{\psi_\alpha}$, then $\zeta - \bE \zeta \in L_{\psi_\alpha}$, and $\|\zeta - \bE \zeta\|_{\psi_\alpha} \lesssim \|\zeta\|_{\psi_\alpha}$. If $\zeta_1, \ldots, \zeta_N$ are $N$ (not necessarily independent nor identically distributed) $\psi_{\alpha}$ random variables defined on $(\Omega, \cA)$, then
\begin{align}\label{eq:product_sub_Weilbull}
    \left\|\prod_{i=1}^N \zeta_i\right\|_{{\psi_\beta}} \leq \prod_{i=1}^N \|\zeta_i\|_{{\psi_{\alpha_i}}},\mbox{ where }\frac{1}{\beta} = \sum_{i=1}^N \frac{1}{\alpha_i}.
\end{align}

The following Lemma is taken from \cite{kuchibhotla_moving_2022}.
\begin{Lemma}\label{lemma:kuchibhotla_sub_Weibull}
    Let $\zeta_1,\cdots,\zeta_N$ be independent mean zero random variables with $\max(\|\zeta_i\|_{\psi_\alpha}:\, i\in[N])<\infty$ for some $\alpha>0$. Then there exist absolute constants $\nC\label{C_sub_Weibull_1}>1$ and $\nC\label{C_sub_Weibull_2}>1$ depending only on $\alpha$ such that for any $t\geq 0$,
    \begin{align*}
        \bP\left(\left|\sum_{i=1}^N\zeta_i\right|\geq \oC{C_sub_Weibull_1}\left(\sum_{i=1}^N \|\zeta_i\|_{\psi_\alpha}^2\right)^{\frac{1}{2}}\sqrt{t} + \oC{C_sub_Weibull_2}t^{\frac{1}{\alpha}}\norm{\vb}_{\beta(\alpha)}\right)\leq 2\exp(-t),
    \end{align*}where $\vb = (\|\zeta_1\|_{\psi_\alpha},\cdots,\|\zeta_N\|_{\psi_\alpha})$ and $\beta(\alpha)=\infty$ when $\alpha\leq 1$ and $\beta(\alpha)=\frac{\alpha}{\alpha-1}$ if $\alpha>1$.
\end{Lemma}
We shall see that the phase transition of $\alpha$ at $1$ is the fundamental reason behind the phase transition at $q = 2$ in the upper bound of the $L^q$ empirical process. Lemma~\ref{lemma:kuchibhotla_sub_Weibull} can be viewed as a generalized Bernstein inequality, which characterizes the tail probability of the sum of independent, mean-zero sub-Weibull random variables. We observe that the tail probability exhibits a phase transition: when $\alpha\leq 2$, the process has a (sub-)Gaussian type tail for small $t$, where $\sqrt{t}$ corresponds to $\exp(-t)$, while for large $t$, it has a  (sub-)Weibull ($\alpha$) type tail, where $t^{1/\alpha}$ corresponds to $\exp(-t)$; see \cite[discussion after Corollary 2.8.3]{vershynin_high-dimensional_2018}. While when $\alpha > 2$, the relation between these two regimes is reversed, that is, it exhibits Gaussian tails for large deviations and Weibull tails for small deviations.

Although, as a sub-Weibull random variable, its concentration inequality is expressed in terms of two metrics—seemingly suggesting that the upper bound of the $L^q$ process should be controlled by the $\gamma$-functionals associated with both norms—the experience (\cite{dirksen_tail_2015,bednorz_concentration_2014,mendelson_upper_2016,al-ghattas_sharp_2025}) with quadratic processes indicates otherwise: if one regards a sub-Weibull random variable as the $q$-th power of a sub-Gaussian random variable, then the sub-Gaussian metric alone suffices to describe the concentration inequality of the $L^q$ process. To see this, applying Lemma~\ref{lemma:kuchibhotla_sub_Weibull} to $\alpha=2/q$ for $1<q<2$, there exists an absolute constant $\nC\label{C_sub_Weibull_3} = \oC{C_sub_Weibull_3}(q) > 1$ such that for any $t\geq 0$, we have
\begin{align*}
    &\bP\left(\left| \sum_{i=1}^N f^q(X_i)-\bE f^q(X) \right| \geq \oC{C_sub_Weibull_1}\left(\sum_{i=1}^N \|f^q(X_i)\|_{\psi_\frac{2}{q}}^2\right)^{\frac{1}{2}}\sqrt{t} + \oC{C_sub_Weibull_2}t^{\frac{q}{2}}(\sum_{i=1}^N \|f^q(X_i)\|_{\psi_{2/q}}^{\frac{2}{2-q}})^{\frac{2-q}{2}} \right) \\
    & \leq \bP\left(\left|\sum_{i=1}^N f^q(X_i)-\bE f^q(X) \right| \geq \oC{C_sub_Weibull_3}\left(\sum_{i=1}^N \|f(X_i)\|_{\psi_2}^{2q}\right)^{\frac{1}{2}}\sqrt{t} + \oC{C_sub_Weibull_3}t^{\frac{q}{2}}(\sum_{i=1}^N \|f(X_i)\|_{\psi_2}^{\frac{2q}{2-q}})^{\frac{2-q}{2}} \right) \leq 2\exp(-t).
\end{align*}Since $X_1,\cdots,X_N$ are i.i.d., there exists an absolute constant $\nC\label{C_sub_Weibull_4} = \oC{C_sub_Weibull_4}(q)>1$ such that, for any $f\in\cF$, and for any $x\geq 0$, we have
\begin{align}\label{eq:Bernstein_one_function_sub_Weibull}
    \bP\left\{\left|\frac{1}{N}\sum_{i=1}^N f^q(X_i)-\bE f^q(X)\right|\ge x\norm{f(X)}_{\psi_2}^q\right\}
\le 2\exp\Bigl(-\frac{1}{\oC{C_sub_Weibull_4}}N\min\Bigl\{x^2,\;x^\frac{2}{q}\Bigr\}\Bigr).
\end{align}
Replacing $f$ with $|f|$, it follows from Fubini's theorem that there exists an absolute constant $\nC\label{C_upper_moment_1}$ depending only on $q$, such that for any $r \ge 1$, one has
\begin{align}\label{eq:upper_moment_1}
    \left(\bE\left|\frac{1}{N}\sum_{i=1}^N |f(X_i)|^q-\bE [|f(X)|^q]\right|^r\right)^{\frac{1}{r}} \leq \oC{C_upper_moment_1}\norm{f}_{\psi_2}^q\left(\sqrt{\frac{r}{N}} + \left(\frac{r}{N}\right)^{\frac{q}{2}}\right).
\end{align}
\begin{Remark}\label{remark:optimality}
    Let $\cF = \{\vzero, \pm f\}$, where $f \in L_{\psi_2}$. Then the $L^q$ process degenerates to the generalized Bernstein inequality \eqref{eq:Bernstein_one_function_sub_Weibull}. Since \cite{kuchibhotla_moving_2022} proved that this inequality is sharp, Theorem~\ref{theo:main} is therefore also sharp.
\end{Remark}

For the $L^q$ empirical process, its tail probability exhibits a phase transition, corresponding to the two types of tail behavior described in \eqref{eq:Bernstein_one_function_sub_Weibull}. Therefore, to ensure that the tail probability given by \eqref{eq:Bernstein_one_function_sub_Weibull} balances the metric entropy of $T_n$, that is, to make the union bound applicable to all elements in $T_{n-1}\times T_n$, we require $|T_{n-1}||T_n| \leq N_{n+1}$ not to exceed $2\exp\!\Bigl(\oC{C_sub_Weibull_4}^{-1}N \min\bigl\{ x^2,x^{2/q} \bigr\}\Bigr)$. The corresponding balance points are: (a) for the Gaussian-type tail, $N_{n+1} \sim \exp(N x^2)$, which is equivalent to $x \sim \sqrt{\log(N_{n+1})/N} \sim (2^n/N)^{1/2}$, denoted by $x^{(G)}_{N,n}$; and (b) for the Weibull-type tail, $N_{n+1} \sim \exp(N x^{2/q})$, which is equivalent to $x \sim (2^n/N)^{q/2}$, denoted by $x^{(W)}_{N,n}$. We need $x$ to be large enough so that both conditions are satisfied, that is, $x\geq \max\{ x^{(G)}_{N,n},\, x^{(W)}_{N,n} \}$. The boundary between these two regimes occurs when $2^n \sim N$. Note that this boundary point is independent of $q$. This boundary is crucial for generic chaining, and we denote by $\fm^*$ the integer such that $2^{\fm^*} \le N < 2^{\fm^*+1}$, i.e., $2^{\fm^*/2} \sim \sqrt{N}$. We refer to this moment $\fm^*$ as the critical time. When $n \le \fm^*$, we have $2^n \le N$, hence $x^{(G)}_{N,n}$ dominates and the tail probability is governed by the Gaussian regime; when $n > \fm^*$, it is dominated by the Weibull regime. In other words, during the chain growth up to time $\fm^*$, the Gaussian tail balances the metric entropy, while beyond $\fm^*$, the Weibull tail takes over.

For any $f \in \cF$, we divide its corresponding chain into two parts: before and after $\fm^*$. That is, $f = (f - \pi_{\fm^*} f) + \pi_{\fm^*} f$. We refer to these two segments of the chain as the terminal part and the initial part, respectively. We now examine, on the two segments of the chain, how the process margin varies for the index increment $\pi_n f - \pi_{n-1} f$.

\subsubsection{Stochastic Argument}

In this section, we study the tail probability of the increment of the margin of the stochastic process $|(P_N - P)|f|^q|$ corresponding to the index increment $\pi_n f - \pi_{n-1} f$. We first examine the initial part of the chain.
\paragraph{Initial part of the chain.}
Applying the Lagrange mean value theorem to the function $x \in \bR \mapsto x^q \in \bR$ over the interval $[a,b]$ where $a,b\geq 0$, we obtain that there exists $\xi \in [a,b]$ (or $\xi \in [b,a]$) such that $a^q - b^q = q\xi^{q-1}(a-b)$ and $|a^q - b^q| \leq q |\xi|^{q-1} |a - b|\leq q\max\{a^{q-1},b^{q-1}\}|a-b|$. Taking the $\psi_{\frac{2}{q}}$ norm on both sides and using \eqref{eq:product_sub_Weilbull} twice together with triangular inequality and identifying $a$, $b$ by $|f|$ and $|g|$ respectively, we obtain
\begin{align*}
    &\norm{|f|^q-|g|^q}_{\psi_{\frac{2}{q}}}\leq q\norm{\max\{|f|^{q-1},|g|^{q-1}\}||f|-|g||}_{\psi_{\frac{2}{q}}}\leq q\norm{\max\{|f|^{q-1},|g|^{q-1}\}}_{\psi_{\frac{2}{q-1}}} \norm{f-g}_{\psi_2} \\
    &\leq q\left(\norm{f}_{\psi_2}^{q-1}+\norm{g}_{\psi_2}^{q-1}\right)d_{\psi_2}(f,g),
\end{align*}and consequently, by Lemma~\ref{lemma:kuchibhotla_sub_Weibull} applied to $\zeta_i = |f(X_i)|^q-|g(X_i)|^q - \bE[|f(X)|^q-|g(X)|^q]$, and $\alpha = \frac{2}{q}$ we obtain that there exists an absolute constant $\nC\label{C_Gaussian}>1$ depending only on $q$ such that for any $x\geq 0$,
\begin{align*}
    \bP\left(\bigg|\frac{1}{N}\sum_{i=1}^N |f(X_i)|^q - |g(X_i)|^q - \bE\left[|f(X)|^q-|g(X)|^q\right]\bigg|\geq \oC{C_Gaussian} d_{\psi_2}(f,g)\diam(\cF)^{q-1}\left(\frac{\sqrt{x}}{\sqrt{N}}+\frac{x^{q/2}}{N^{q/2}}\right)\right)\leq 2\exp(-x).
\end{align*}Letting $x=2^nu$ with $u\geq 1$ and $n<\fm^*$, then $(2^n/N)^{q/2}\leq (2^n/N)^{1/2}$. Therefore, there exists an absolute constant $\nC\label{C_Gaussian_1}>1$ depending only on $q$ such that
\begin{align}\label{eq:increment_sub_Gaussian}
    \forall n<\fm^*,\, \forall u\geq 1,\, \forall f,g\in L_{\psi_2},\, \bP\left(\bigg|(P-P_N)(|f|^q - |g|^q)\bigg|\geq \oC{C_Gaussian_1}u^{q/2} 2^{n/2}d_{\psi_2}(f,g)\frac{\diam(\cF)^{q-1}}{\sqrt{N}}\right)\leq 2\exp(-u2^n).
\end{align}Here, when $n \ge \fm^*$, the term $(2^n / N)^{q/2}$ dominates. Hence, if we set $x = 2^n u$ with $u \ge 1$ and $n \ge \fm^*$, then the high-probability upper bound for the increment of the process margin $(P_N - P)|f|^q$ contains a factor of $2^{nq/2}$. This would lead to the appearance of the $\gamma_{2/q}(\cF, \psi_2)$ functional in the final generic chaining bound, which we aim to avoid (since we want the final result to involve only the $\gamma_2$ functional). Therefore, for $n \ge \fm^*$, we need to adopt a different treatment.

\paragraph{Terminal part of the chain.}
Replacing $f$ in equation~\eqref{eq:Bernstein_one_function_sub_Weibull} with $f - g$, where $f, g \in L_{\psi_2}$, and adjusting the value of $x$, we note that by the sub-Gaussian increment assumption, there exists some absolute constant $C(q)\leq \frac{\sqrt{\pi}}{2}\oC{C_subgaussian}\sqrt{q}$, such that for any $f,g\in L_{\psi_2}$, there holds $\|f - g\|_{L^q} \leq C(q)d_{\psi_2}(f, g)$, (for example, this can be obtained from the following Lemma~\ref{lemma:dirksen}). Hence, there exists an absolute constant $\nC\label{C_Weibull} = \oC{C_Weibull}(q) > 1$ such that for any $f, g \in L_{\psi_2}$ and any $x \ge 1$, one has
\begin{align*}
    \bP\left\{ \left(\frac{1}{N}\sum_{i=1}^N |f(X_i)-g(X_i)|^q\right)^{\frac{1}{q}}\geq\oC{C_Weibull} x\norm{f-g}_{\psi_2}\right\}
\le 2\exp\Bigl(-Nx^2\Bigr).
\end{align*}Let $x = \sqrt{\tfrac{u}{N}}\, 2^{n/2}$, where $n \ge \fm^*$ and $u \ge 1$, then
\begin{align}\label{eq:increment_sub_Weibull}
    \forall n\geq\fm^*,\, \forall u\geq 1,\,\forall f,g\in L_{\psi_2},\, \bP\left\{ \left(\frac{1}{N}\sum_{i=1}^N |f(X_i)-g(X_i)|^q\right)^{\frac{1}{q}}\geq\oC{C_Weibull} \sqrt{\frac{u}{N}}2^{n/2}\norm{f-g}_{\psi_2}\right\}
\le 2\exp\Bigl(-u2^n\Bigr).
\end{align}We also need the following lemma, taken from \cite[Lemma A.5]{dirksen_tail_2015}.
\begin{Lemma}[\cite{dirksen_tail_2015}]\label{lemma:dirksen}
    Fix $1\leq q<\infty$ and $0<\alpha<\infty$. Let $\gamma\geq0$ and suppose that $\xi$ is a positive random variable such that for some $c\geq 1$ and $t_*>0$,
    \begin{align*}
        \forall t\geq t_*,\, \bP\big(\xi>\gamma t\big)\leq c\exp\big(-\frac{1}{4}qt^\alpha\big).
    \end{align*}Then there exists an absolute constant $c_\alpha>0$ depending only on $\alpha$ such that $\|\xi\|_{L^q}\leq \gamma(c_\alpha c + t_*)$.
\end{Lemma}
We now apply Lemma~\ref{lemma:dirksen} to $\xi = \bigl(\tfrac{1}{N}\sum_{i=1}^N |f(X_i) - g(X_i)|^q\bigr)^{1/q}$. Let $u = (q/4) 2^{-n} t^2$ in equation~\eqref{eq:increment_sub_Weibull}, where $t$ is as in Lemma~\ref{lemma:dirksen}. Then, from equation~\eqref{eq:increment_sub_Weibull}, we have $\alpha = 2$, $c = 2$, $\gamma = \oC{C_Weibull}\sqrt{q}/(2\sqrt{N})\, d_{\psi_2}(f,g)$, and $t_* = 2^{1 + n/2} q^{-1/2}$, which satisfy the conditions of Lemma~\ref{lemma:dirksen}. Consequently, Lemma~\ref{lemma:dirksen} yields
$\|\xi\|_{L^q} \le \tfrac{\oC{C_Weibull}}{\sqrt{N}}\, d_{\psi_2}(f,g)\, (2^{n/2} + c_2\sqrt{q})$. That is, there exists an absolute constant $\nC\label{C_Weibull_upper_moment} = \oC{C_Weibull_upper_moment}(q)>1$ such that
\begin{align}\label{eq:increment_sub_Weibull_moment}
    &\forall n\geq \fm^*,\, \forall f,g\in L_{\psi_2},\, \quad \left(\bE\left[\frac{1}{N}\sum_{i=1}^N\big|f(X_i)-g(X_i)\big|^q\right]\right)^{\frac{1}{q}}\leq \frac{\oC{C_Weibull_upper_moment}}{\sqrt{N}}\, 2^{n/2}d_{\psi_2}(f,g).
\end{align}
Here, we do not directly analyze the increments of the margin of the $L^q$ process, but this is already sufficient for applying generic chaining. We employ Dirksen's contraction technique \cite{dirksen_tail_2015} in this step.

\paragraph{Summary of stochastic arguments.}Following the terminology of \cite{dirksen_tail_2015}, we define the random event $\Omega_{u,p}$ as
\begin{align}\label{eq:def_Omega_u,p}
        \Omega_{u,p} := \Bigg\{ \forall f\in\cF,\, \begin{cases}
            \bigg|(P-P_N)(|\pi_nf|^q - |\pi_{n-1}f|^q)\bigg|\leq \oC{C_Gaussian_1}u^{q/2} 2^{n/2}d_{\psi_2}(\pi_nf,\pi_{n-1}f)\frac{\diam(\cF)^{q-1}}{\sqrt{N}},\, &\forall \ell<n<\fm^*,\\
            \left(\frac{1}{N}\sum_{i=1}^N |\pi_nf(X_i)-\pi_{n-1}f(X_i)|^q\right)^{\frac{1}{q}}\leq\oC{C_Weibull} \sqrt{\frac{u}{N}}2^{n/2}\norm{\pi_nf-\pi_{n-1}f}_{\psi_2}, & \forall n\geq \fm^*
        \end{cases}\Bigg\},
\end{align}where we recall that $\ell=\lfloor \log_2(p)\rfloor$. 
The following lemma is taken from \cite[Lemma A.4]{dirksen_tail_2015}.
\begin{Lemma}[\cite{dirksen_tail_2015}]\label{lemma:union_bound}
    Fix $1\leq p<\infty$, $0<\alpha<\infty$, $u\geq 2^{1/\alpha}$ and set $\ell=\lfloor \log_2(p)\rfloor$. For every $n>\ell$, let $(\Omega_i^{(n)})_{i\in I_n}$ be a collection of events satisfying
    \begin{align*}
        \bP\left(\Omega_i^{(n)}\right)\leq 2\exp(-2^n u^\alpha),\quad \forall i\in I_n.
    \end{align*}If $|I_n|\leq N_{n+1}$, then there exists an absolute constant $\nC\label{C_prob}\leq 17$ such that
    \begin{align*}
        \bP\left( \bigcup_{n>\ell}\bigcup_{i\in I_n}\Omega_i^{(n)} \right)\leq \oC{C_prob}\exp\left(-\frac{1}{4}pu^\alpha\right).
    \end{align*}
\end{Lemma}
By Lemma~\ref{lemma:union_bound}, we obtain that for any $u\geq 2$ and $1\leq p<\infty$, there exists an absolute constant $\oC{C_prob}\leq 17$ such that
\begin{align}\label{eq:prob_Omega}
    \bP(\Omega_{u,p}) \geq 1-\oC{C_prob}\exp(-pu/4).
\end{align}

\subsubsection{Deterministic Argument}
In this section, we work on the event $\Omega_{u,p}$.

Define $I_{\mathrm{Gauss}} = \{\ell+1,\cdots,\fm^*\}$, and $I_{\mathrm{Weibull}} = \{\fm^*+1,\cdots\}$. Then for any $f\in\cF$,
\begin{align}\label{eq:chain_expansion}
    (P - P_N)|f|^q = (P-P_N)|\pi_\ell f|^q + \sum_{n\in I_{\mathrm{Gauss}}}(P-P_N)(|\pi_n f|^q - |\pi_{n-1}f|^q) + \sum_{n\in I_{\mathrm{Weibull}}}(P-P_N)(|\pi_n f|^q - |\pi_{n-1}f|^q).
\end{align}We first handle the terminal part of the chain. When $q > 1$, this corresponds to the Weibull component.

\subsubsection{The terminal of the chain, \texorpdfstring{$I_{\mathrm{Weibull}}$}{Weibull}} We use triangular inequality
\begin{align}\label{eq:Weibull_chaining}
    \begin{aligned}
        &\left|\sum_{n\in I_{\mathrm{Weibull}}}(P-P_N)\left(|\pi_n f|^q - |\pi_{n-1}f|^q\right)\right|\\
        &\leq \sum_{n\in I_{\mathrm{Weibull}}}\left|P_N\left(|\pi_n f|^q - |\pi_{n-1}f|^q\right)\right|+\sum_{n\in I_{\mathrm{Weibull}}}\left|P\left(|\pi_n f|^q - |\pi_{n-1}f|^q\right)\right|.
    \end{aligned}
\end{align}In the following two parts, we separately handle the upper bounds for the empirical and population components.

\paragraph{Upper bound for the empirical part.}
Applying the Lagrange mean value theorem to the function $x \in \bR \mapsto x^q \in \bR$ over the interval $[a,b]$, we obtain that there exists $\xi \in [a,b]$ (or $\xi \in [b,a]$) such that $a^q - b^q = q\xi^{q-1}(a-b)$. Applying to $a=|(\pi_nf)(X_i)|$ and $b = |(\pi_{n-1}f)(X_i)|$ for each $1\leq i\leq N$, we obtain that $|(\pi_nf)(X_i)|^q - |(\pi_{n-1}f)(X_i)|^q = q\xi_i^{q-1}(|(\pi_nf)(X_i)| - |(\pi_{n-1}f)(X_i)|)$, where $\xi_i\in[|(\pi_nf)(X_i)|,|(\pi_{n-1}f)(X_i)|]$ (or vice versa). Taking sum over $1\leq i\leq N$ and dividing by $N$, we have $P_N(|\pi_nf|^q - |\pi_{n-1}f|^q) = qP_N\xi^{q-1}(|\pi_nf| - |\pi_{n-1}f|)$. Taking absolute value, using H{\"o}lder inequality for $(q,\frac{q}{q-1})$ and by using that $|\xi_i|\leq \max\{|(\pi_nf)(X_i)|,|(\pi_{n-1}f)(X_i)|\}$, we have
\begin{align}\label{eq:empirical_Weibull_upper}
    \begin{aligned}
        &\left|P_N \left(|\pi_n f|^q - |\pi_{n-1}f|^q\right)\right| = q\left|\frac{1}{N}\sum_{i=1}^N \xi_i^{q-1}\left(|(\pi_n f)(X_i)| - |(\pi_{n-1}f)(X_i)|\right)\right|\\
        &\leq q\left(\frac{1}{N}\sum_{i=1}^N\left|\xi_i\right|^q\right)^{\frac{q-1}{q}}\left( \frac{1}{N}\sum_{i=1}^N \bigg|\left|(\pi_n f)(X_i)\right| - \left|(\pi_{n-1}f)(X_i)\right|\bigg|^q\right)^{\frac{1}{q}}  \\
        &\leq    q\left(\frac{1}{N}\sum_{i=1}^N |(\pi_nf)(X_i)|^q+|(\pi_{n-1}f)(X_i)|^q \right)^{\frac{q-1}{q}}\left( \frac{1}{N}\sum_{i=1}^N \bigg|\left|(\pi_n f)(X_i)\right| - \left|(\pi_{n-1}f)(X_i)\right|\bigg|^q\right)^{\frac{1}{q}}.
    \end{aligned}
\end{align}Now, we recognize that
\begin{align*}
    &\frac{1}{N}\sum_{i=1}^N |(\pi_nf)(X_i)|^q - \bE|(\pi_nf)(X)|^q = (P_N-P)|\pi_nf|^q.
\end{align*}Therefore, by the elementary inequality $(a+b)^{\frac{1}{q}}\leq  (a^{\frac{1}{q}}+b^{\frac{1}{q}})$,
\begin{align*}
    &\left(\frac{1}{N}\sum_{i=1}^N |(\pi_nf)(X_i)|^q+|(\pi_{n-1}f)(X_i)|^q \right)^{\frac{1}{q}}\\ 
    &= \left( \left[(P_N-P)|\pi_nf|^q + P|\pi_nf|^q\right] + \left[(P_N-P)|\pi_{n-1}f|^q + P|\pi_{n-1}f|^q\right] \right)^{\frac{1}{q}}\\
    &\leq  \left( \left| P_N-P\right||\pi_nf|^q + P|\pi_nf|^q\right)^{\frac{1}{q}} +  \left(  \left|P_N-P\right||\pi_{n-1}f|^q + P|\pi_{n-1}f|^q \right)^{\frac{1}{q}}\\
    &\leq 2\left( \sup_{f\in\cF}\bigg|(P-P_N)|f|^q\bigg| + \diam(\cF)^q \right)^{\frac{1}{q}}.
\end{align*}Taking $q-1$ power, then
\begin{align*}
    &q\left(\frac{1}{N}\sum_{i=1}^N |(\pi_nf)(X_i)|^q+|(\pi_{n-1}f)(X_i)|^q \right)^{\frac{q-1}{q}}\lesssim_q \sup_{f\in\cF}\bigg|\big(P-P_N\big)|f|^q\bigg|^{\frac{q-1}{q}} + \diam(\cF)^{q-1}.
\end{align*}Together with \eqref{eq:def_Omega_u,p}, we obtain
\begin{align}\label{eq:Weibull_P_N_result}
    &\bigg|P_N \left(|\pi_n f|^q - |\pi_{n-1}f|^q\right)\bigg| \lesssim_q \oC{C_Weibull} \sqrt{u}\frac{2^{\frac{n}{2}}}{\sqrt{N}}\norm{f-g}_{\psi_2}\left(\sup_{f\in\cF}\left|\big(P-P_N\big)|f|^q\right|^{\frac{q-1}{q}} + \diam(\cF)^{q-1}\right).
\end{align}
\paragraph{Upper bound for the population part.}Recall that $\{X_i\}_{i=1}^N$ are i.i.d., and hence the population part can be equivalently written as
\begin{align*}
    &\bigg|P\left(|\pi_n f|^q - |\pi_{n-1}f|^q\right)\bigg| = \left| \bE\left[\frac{1}{N}\sum_{i=1}^N |\pi_nf(X_i)|^q - |\pi_{n-1}f(X_i)|^q \right]\right|\leq  \bE\left|\frac{1}{N}\sum_{i=1}^N |\pi_nf(X_i)|^q - |\pi_{n-1}f(X_i)|^q \right|,
\end{align*}where we have used Jensen's inequality. By \eqref{eq:empirical_Weibull_upper},
\begin{align*}
    &\bE\left|\frac{1}{N}\sum_{i=1}^N |\pi_nf(X_i)|^q - |\pi_{n-1}f(X_i)|^q \right|\\
    &\leq q\bE\left[\left(\frac{1}{N}\sum_{i=1}^N |(\pi_nf)(X_i)|^q+|(\pi_{n-1}f)(X_i)|^q \right)^{\frac{q-1}{q}}\left( \frac{1}{N}\sum_{i=1}^N \bigg|\left|(\pi_n f)(X_i)\right| - \left|(\pi_{n-1}f)(X_i)\right|\bigg|^q\right)^{\frac{1}{q}}\right]\\
    &\leq q\left(\bE\left[\frac{1}{N}\sum_{i=1}^N |(\pi_nf)(X_i)|^q+|(\pi_{n-1}f)(X_i)|^q \right]\right)^{\frac{q-1}{q}}\left(\bE\left[\frac{1}{N}\sum_{i=1}^N \bigg|\left|(\pi_n f)(X_i)\right| - \left|(\pi_{n-1}f)(X_i)\right|\bigg|^q\right]\right)^{\frac{1}{q}}\\
    &\leq q2^{\frac{q-1}{q}}\diam(\cF)^{q-1} \left(\bE\left[\frac{1}{N}\sum_{i=1}^N \bigg|\left|(\pi_n f)(X_i)\right| - \left|(\pi_{n-1}f)(X_i)\right|\bigg|^q\right]\right)^{\frac{1}{q}},
\end{align*}where we have used H{\"o}lder inequality with conjugate pair $(q,\frac{q}{q-1})$. Applying \eqref{eq:increment_sub_Weibull_moment}, we obtain that
\begin{align*}
    &\left(\bE\left[\frac{1}{N}\sum_{i=1}^N \bigg|\left|(\pi_n f)(X_i)\right| - \left|(\pi_{n-1}f)(X_i)\right|\bigg|^q\right]\right)^{\frac{1}{q}}\leq \oC{C_Weibull_upper_moment} 2^{n/2}\|\pi_nf-\pi_{n-1}f\|_{\psi_2}\frac{1}{\sqrt{N}}.
\end{align*}As a result,
\begin{align}\label{eq:population_Weibull_upper}
    &\bigg|P\left(|\pi_n f|^q - |\pi_{n-1}f|^q\right)\bigg|\lesssim_q \frac{\diam(\cF)^{q-1}}{\sqrt{N}}2^{n/2}d_{\psi_2}(\pi_nf,\pi_{n-1}f).
\end{align}

\paragraph{Upper bound for the terminal of the chain.} Combining \eqref{eq:Weibull_chaining}, \eqref{eq:Weibull_P_N_result}, and \eqref{eq:population_Weibull_upper}, we obtain the upper bound for the terminal part of the chain. That is, on the event $\Omega_{u,p}$, there exists an absolute constant $\nC\label{C_upper_Weibull_result} > 1$ depending only on $q$, such that for any $f \in \cF$, one has
\begin{align}\label{eq:result_terminal_chain}
    \begin{aligned}
        &\left|\sum_{n\in I_{\mathrm{Weibull}}}(P-P_N)\left(|\pi_n f|^q - |\pi_{n-1}f|^q\right)\right| \leq \oC{C_upper_Weibull_result} \Bigg( \sup_{f\in\cF}\left|\big(P-P_N\big)|f|^q\right|^{\frac{q-1}{q}} + \diam(\cF)^{q-1}\Bigg) \frac{\sqrt{u}}{\sqrt{N}}\gamma_2(\cF,d_{\psi_2}),
    \end{aligned}
\end{align}where we have used that $\sum_{n\in I_{\mathrm{Weibull}}}2^{n/2}d_{\psi_2}(\pi_nf,\pi_{n-1}f)\leq \gamma_2(\cF,d_{\psi_2})$.

\begin{Remark}
    An interesting observation is that, when dealing with the terminal segment of the chain, we in fact only need to use $\diam(\cF, L^q)$, rather than $\diam(\cF, L_{\psi_2})$. We keep the notation $\diam(\cF) = \diam(\cF, L_{\psi_2})$ because, in the initial segment below, we will need $\diam(\cF, L_{\psi_2})$. If one could replace $\diam(\cF) = \diam(\cF, L_{\psi_2})$ in \eqref{eq:increment_sub_Gaussian} by $\diam(\cF, L^q)$, then the resulting upper bound would depend only on $\diam(\cF, L^q)$ and not on $\diam(\cF, L_{\psi_2})$ (this would be analogous to the result of \cite[Theorem 1.13]{mendelson_upper_2016}). However, at present we do not know how to perform such a replacement.
\end{Remark}

\subsubsection{The initial of the chain, \texorpdfstring{$I_{\mathrm{Gauss}}$}{Gauss}}

For the initial part of the chain, we directly use the increments of the process margins provided by $\Omega_{u,p}$ for $\ell < n < \fm^*$, together with the triangle inequality, to obtain
\begin{align}\label{eq:Gaussian_chain_result}
    &\left|\sum_{n\in I_{\mathrm{Gauss}}}(P-P_N)|\pi_n f|^q - |\pi_{n-1}f|^q\right| \leq \oC{C_Gaussian_1}u^{q/2}\frac{\diam(\cF)^{q-1}}{\sqrt{N}}\gamma_2(\cF,d_{\psi_2}),
\end{align}where we have used that $\sum_{n\in I_{\mathrm{Gauss}}}2^{n/2}d_{\psi_2}(\pi_nf,\pi_{n-1}f)\leq \gamma_2(\cF,d_{\psi_2})$.

\subsubsection{Combining two segments}

Combining \eqref{eq:chain_expansion}, \eqref{eq:result_terminal_chain} and \eqref{eq:Gaussian_chain_result}, we obtain that on $\Omega_{u,p}$, for any $f\in\cF$, one has
\begin{align*}
    &\sup_{f\in\cF}\bigg|(P-P_N)|f|^q\bigg|\leq \sup_{f\in\cF}\bigg|(P-P_N)|\pi_\ell f|^q\bigg|+  (\oC{C_Gaussian_1}+\oC{C_upper_Weibull_result})u^{q/2}\diam(\cF)^{q-1}\frac{\gamma_2(\cF,d_{\psi_2})}{\sqrt{N}}\\
    &+\oC{C_upper_Weibull_result}\sqrt{u}\frac{\gamma_2(\cF,d_{\psi_2})}{\sqrt{N}}\left(\sup_{f\in\cF}\bigg|\big(P-P_N\big)|f|^q\bigg|\right)^{\frac{q-1}{q}},
\end{align*}where we have used the fact that $u^{1/2}\leq u^{q/2}$ when $q\geq 1$ and $u\geq 1$.
This yields an inequality for $\sup_{f \in \cF} \big| (P - P_N)|f|^q \big|^{1/q}$. We use the following elementary lemma:
\begin{Lemma}\label{lemma:elementary}
    Suppose $Y,a,b,S\geq 0$, $q\geq 1$. If $Y^q\leq aY^{q-1} + (b+S)$. Then $Y - b^{1/q}\leq a + S^{1/q}$.
\end{Lemma}
\beginproof
We first prove that $Y\leq a+(b+S)^{1/q}$. Suppose $Y>a+(b+S)^{1/q}$, then $Y>(b+S)^{1/q}$ and $Y^{q-1}>(b+S)^{(q-1)/q}$. Therefore $Y^{q-1}(Y-a)>b+S$. However, from our assumption, $Y^q\leq aY^{q-1} + (b+S)$, that is, $Y^{q-1}(Y-a)\leq b+S$, which is a contradiction. Hence we must have $Y\leq a+(b+S)^{1/q}$. Now since $b,S\geq 0$ and $q\geq 1$, we use the elementary inequality $(b+S)^{1/q}\leq b^{1/q} + S^{1/q}$ to conclude the proof.
\endproof
In Lemma~\ref{lemma:elementary}, let $Y = \sup_{f \in \cF} \big| (P - P_N)|f|^q \big|^{1/q}$, $a = \oC{C_upper_Weibull_result}\sqrt{u}\,\tfrac{\gamma_2(\cF,d_{\psi_2})}{\sqrt{N}}$, $b = \sup_{f\in\cF}\big|(P-P_N)|\pi_\ell f|^q\big|$, and $S = (\oC{C_Gaussian_1}+\oC{C_upper_Weibull_result})u^{q/2}\,\diam(\cF)^{q-1}\tfrac{\gamma_2(\cF,d_{\psi_2})}{\sqrt{N}}$, we obtain that on $\Omega_{u,p}$, for any $f\in\cF$,
\begin{align*}
    &\sup_{f \in \cF} \big| (P - P_N)|f|^q \big|^{1/q} - \sup_{f\in\cF}\big|(P-P_N)|\pi_\ell f|^q\big|^{1/q} \\
    &\leq \sqrt{u}\Bigg(\oC{C_upper_Weibull_result}\,\frac{\gamma_2(\cF,d_{\psi_2})}{\sqrt{N}} + (\oC{C_Gaussian_1}+\oC{C_upper_Weibull_result})^{1/q}\,\diam(\cF)^{\frac{q-1}{q}}\frac{\gamma_2^{1/q}(\cF,d_{\psi_2})}{N^{1/(2q)}}\Bigg).
\end{align*}
Since $\pi_\ell f \in \cF$, the left-hand side $\sup_{f \in \cF} \big|(P - P_N)|f|^q\big|^{1/q} - \sup_{f \in \cF} \big|(P - P_N)|\pi_\ell f|^q\big|^{1/q}$ is nonnegative. Moreover, by \eqref{eq:prob_Omega}, we may apply Lemma~\ref{lemma:dirksen} with $\xi$ to be the left-hand-side of the above inequality, $q = p$, $t = \sqrt{u}$ (hence $t_* = \sqrt{2}$), $\alpha = 2$, $c = \oC{C_prob}$, and $\gamma$ equal to the right-hand side above divided by $\sqrt{u}$, to obtain
\begin{align*}
    &\Bigg(\bE\left[\bigg(\sup_{f \in \cF} \big| (P - P_N)|f|^q \big|^{1/q} - \sup_{f\in\cF}\big|(P-P_N)|\pi_\ell f|^q\big|^{1/q}\bigg)^p\right]\Bigg)^{\frac{1}{p}} \leq \oC{C_upper_moment}\Bigg(\frac{\gamma_2(\cF,d_{\psi_2})}{\sqrt{N}} + \diam(\cF)^{\frac{q-1}{q}}\frac{\gamma_2^{1/q}(\cF,d_{\psi_2})}{N^{1/(2q)}}\Bigg),
\end{align*}where $\nC\label{C_upper_moment}$ is an absolute constant that depend only on $\oC{C_prob}$, $\oC{C_Gaussian_1}$ and $\oC{C_upper_Weibull_result}$. By the triangle inequality,
\begin{align}\label{eq:pre_concoulsion_1}
    \begin{aligned}
        &\left(\bE\left[\sup_{f \in \cF} \bigg| (P - P_N)|f|^q \bigg|^{\frac{p}{q}}\right]\right)^{\frac{1}{p}}\\
        &\leq \oC{C_upper_moment}\Bigg(\frac{\gamma_2(\cF,d_{\psi_2})}{\sqrt{N}} + \diam(\cF)^{\frac{q-1}{q}}\frac{\gamma_2^{1/q}(\cF,d_{\psi_2})}{N^{1/(2q)}}\Bigg) + \left(\bE\left[\sup_{f\in\cF}\bigg|(P-P_N)|\pi_\ell f|^q\bigg|^{\frac{p}{q}}\right]\right)^{\frac{1}{p}}.
    \end{aligned}
\end{align}The following lemma is taken from \cite[Lemma A.3]{dirksen_tail_2015}.
\begin{Lemma}[\cite{dirksen_tail_2015}]\label{lemma:dirksen_diameter}
    Fix $1\leq p<\infty$, set $\ell=\lfloor\log_2(p)\rfloor$ and let $(X_t)_{t\in T}$ be a collection of complex-valued random variables. If $|T|\leq N_\ell$, then $(\bE\sup_{t\in T}|X_t|^p)^{1/p}\leq 2\sup_{t\in T}(\bE|X_t|^p)^{1/p}$.
\end{Lemma}Applying Lemma~\ref{lemma:dirksen_diameter} to $X_t = ((P_N - P)|f|^q)^{1/q}$ with $T = \cF_\ell$ together with \eqref{eq:upper_moment_1} where $r$ is set to be $p/q$, we obtain that there exists an absolute constant $\nC\label{C_upper_moment_2}$ depending only on $q$ such that
\begin{align*}
    &\left(\bE\left[\sup_{f\in\cF}\bigg|(P-P_N)|\pi_\ell f|^q\bigg|^{\frac{p}{q}}\right]\right)^{\frac{q}{p}}\leq 2^q\sup_{f\in\cF}\left(\bE\left| (P-P_N)|f|^q\right|^{\frac{p}{q}}\right)^{\frac{q}{p}} \leq  \oC{C_upper_moment_2}\diam(\cF)^q\left(\sqrt{\frac{p}{N}} + \left(\frac{p}{N}\right)^{q/2}\right).
\end{align*}
Taking $q$-th power in \eqref{eq:pre_concoulsion_1} on both sides and substituting this into \eqref{eq:pre_concoulsion_1} and rearranging yields the $p$-th moment upper bound \eqref{eq:objective_moments} stated in Theorem~\ref{theo:1<q<2}. Here, we set $\oC{C_final}\geq \max\{\oC{C_upper_moment_2},2^q\oC{C_upper_moment}^q\}$. For the tail probability, we use the following lemma.
\begin{Lemma}\label{lemma:from_moment_to_tail}
    Let $q\geq 1$. Suppose $\zeta$ is a random variable satisfies the following assumption: there exist some parameters $a_0,a_1,a_2\geq 0$ such that for any $p\geq 1$, $\|\zeta\|_{L^p}\leq a_0 + a_1\sqrt{p} + a_2 p^{1/q}$. There then exists an absolute constant $\oC{C_from_moments_to_tail}$ depending only on $q$ such that for any $x\geq 1$,
    \begin{align*}
        \bP\left( |\zeta|\geq e\oC{C_from_moments_to_tail}\left( a_0 + a_1\sqrt{x} + a_2 x^{\frac{1}{q}} \right) \right)\leq \exp(-x).
    \end{align*}
\end{Lemma}
\beginproof
Let $p=\lceil x\rceil$. By Markov's inequality, for any $\lambda>0$, $\bP(|\zeta|\geq \lambda) = \bP(|\zeta|^p\geq \lambda^p)\leq (\frac{\|\zeta\|_{L^p}}{\lambda})^p$. Now, $p<x+1\leq 2x$, hence by assumption, $\|\zeta\|_{L^p}\leq a_0+a_1\sqrt{2x}+a_22^{1/q}x^{1/q}\leq \oC{C_from_moments_to_tail}(a_0+a_1\sqrt{x}+a_2x^{1/q})$, where $\nC\label{C_from_moments_to_tail}=\max\{\sqrt{2},2^{1/q}\}$. Let $\lambda = e\oC{C_from_moments_to_tail}(a_0+a_1\sqrt{x}+a_2x^{1/q})$, then
\begin{align*}
    \bP\left(|\zeta|\geq e\oC{C_from_moments_to_tail}(a_0+a_1\sqrt{x}+a_2x^{1/q})\right)\leq \exp(-p)\leq \exp(-x).
\end{align*}
\endproof
The tail probability \eqref{eq:objective_tail_prob} follows readily by applying Lemma~\ref{lemma:from_moment_to_tail} together with a large enough $\oC{C_final}$. This completes the proof of Theorem~\ref{theo:1<q<2}.


\bibliographystyle{alpha}
\bibliography{biblio}
\end{document}